
\hsize=13.2cm
\vsize=19.5cm
\parindent=0cm   \parskip=0pt
\pageno=1

\def\ind{\hskip 1cm\relax}
\hoffset=15mm
\voffset=1cm

\ifnum\mag=\magstep1
\hoffset=0.2cm   
\voffset=-.5cm   
\fi

\pretolerance=500 \tolerance=1000  \brokenpenalty=5000

\catcode`\@=11

\font\eightrm=cmr8         \font\eighti=cmmi8
\font\eightsy=cmsy8        \font\eightbf=cmbx8
\font\eighttt=cmtt8        \font\eightit=cmti8
\font\eightsl=cmsl8        \font\sixrm=cmr6
\font\sixi=cmmi6           \font\sixsy=cmsy6
\font\sixbf=cmbx6

\skewchar\eighti='177 \skewchar\sixi='177
\skewchar\eightsy='60 \skewchar\sixsy='60

\newfam\gothfam           \newfam\bboardfam
\newfam\cyrfam

\def\tenpoint{%
  \textfont0=\tenrm \scriptfont0=\sevenrm \scriptscriptfont0=\fiverm
  \def\rm{\fam\z@\tenrm}%
  \textfont1=\teni  \scriptfont1=\seveni  \scriptscriptfont1=\fivei
  \def\oldstyle{\fam\@ne\teni}\let\old=\oldstyle
  \textfont2=\tensy \scriptfont2=\sevensy \scriptscriptfont2=\fivesy
  \textfont\gothfam=\tengoth \scriptfont\gothfam=\sevengoth
  \scriptscriptfont\gothfam=\fivegoth
  \def\goth{\fam\gothfam\tengoth}%
  \textfont\bboardfam=\tenbboard \scriptfont\bboardfam=\sevenbboard
  \scriptscriptfont\bboardfam=\sevenbboard
  \def\bb{\fam\bboardfam\tenbboard}%
 \textfont\cyrfam=\tencyr \scriptfont\cyrfam=\sevencyr
  \scriptscriptfont\cyrfam=\sixcyr
  \def\cyr{\fam\cyrfam\tencyr}%
  \textfont\itfam=\tenit
  \def\it{\fam\itfam\tenit}%
  \textfont\slfam=\tensl
  \def\sl{\fam\slfam\tensl}%
  \textfont\bffam=\tenbf \scriptfont\bffam=\sevenbf
  \scriptscriptfont\bffam=\fivebf
  \def\bf{\fam\bffam\tenbf}%
  \textfont\ttfam=\tentt
  \def\tt{\fam\ttfam\tentt}%
  \abovedisplayskip=12pt plus 3pt minus 9pt
  \belowdisplayskip=\abovedisplayskip
  \abovedisplayshortskip=0pt plus 3pt
  \belowdisplayshortskip=4pt plus 3pt
  \smallskipamount=3pt plus 1pt minus 1pt
  \medskipamount=6pt plus 2pt minus 2pt
  \bigskipamount=12pt plus 4pt minus 4pt
  \normalbaselineskip=12pt
  \setbox\strutbox=\hbox{\vrule height8.5pt depth3.5pt width0pt}%
  \let\bigf@nt=\tenrm       \let\smallf@nt=\sevenrm
  \normalbaselines\rm}

\def\eightpoint{%
  \textfont0=\eightrm \scriptfont0=\sixrm \scriptscriptfont0=\fiverm
  \def\rm{\fam\z@\eightrm}%
  \textfont1=\eighti  \scriptfont1=\sixi  \scriptscriptfont1=\fivei
  \def\oldstyle{\fam\@ne\eighti}\let\old=\oldstyle
  \textfont2=\eightsy \scriptfont2=\sixsy \scriptscriptfont2=\fivesy
  \textfont\gothfam=\eightgoth \scriptfont\gothfam=\sixgoth
  \scriptscriptfont\gothfam=\fivegoth
  \def\goth{\fam\gothfam\eightgoth}%
  \textfont\cyrfam=\eightcyr \scriptfont\cyrfam=\sixcyr
  \scriptscriptfont\cyrfam=\sixcyr
  \def\cyr{\fam\cyrfam\eightcyr}%
  \textfont\bboardfam=\eightbboard \scriptfont\bboardfam=\sevenbboard
  \scriptscriptfont\bboardfam=\sevenbboard
  \def\bb{\fam\bboardfam}%
  \textfont\itfam=\eightit
  \def\it{\fam\itfam\eightit}%
  \textfont\slfam=\eightsl
  \def\sl{\fam\slfam\eightsl}%
  \textfont\bffam=\eightbf \scriptfont\bffam=\sixbf
  \scriptscriptfont\bffam=\fivebf
  \def\bf{\fam\bffam\eightbf}%
  \textfont\ttfam=\eighttt
  \def\tt{\fam\ttfam\eighttt}%
  \abovedisplayskip=9pt plus 3pt minus 9pt
  \belowdisplayskip=\abovedisplayskip
  \abovedisplayshortskip=0pt plus 3pt
  \belowdisplayshortskip=3pt plus 3pt
  \smallskipamount=2pt plus 1pt minus 1pt
  \medskipamount=4pt plus 2pt minus 1pt
  \bigskipamount=9pt plus 3pt minus 3pt
  \normalbaselineskip=9pt
  \setbox\strutbox=\hbox{\vrule height7pt depth2pt width0pt}%
  \let\bigf@nt=\eightrm     \let\smallf@nt=\sixrm
  \normalbaselines\rm}

\def\pc#1{\bigf@nt#1\smallf@nt}         \def\pd#1 {{\pc#1} }

\def\^#1{\if#1i{\accent"5E\i}\else{\accent"5E #1}\fi}
\def\"#1{\if#1i{\accent"7F\i}\else{\accent"7F #1}\fi}

\newtoks\auteurcourant      \auteurcourant={\hfil}
\newtoks\titrecourant       \titrecourant={\hfil}

\newtoks\hautpagetitre      \hautpagetitre={\hfil}
\newtoks\baspagetitre       \baspagetitre={\hfil}

\newtoks\hautpagegauche
\hautpagegauche={\eightpoint\rlap{\folio}\hfil\the\auteurcourant\hfil}
\newtoks\hautpagedroite
\hautpagedroite={\eightpoint\hfil\the\titrecourant\hfil\llap{\folio}}

\newtoks\baspagegauche      \baspagegauche={\hfil}
\newtoks\baspagedroite      \baspagedroite={\hfil}

\newif\ifpagetitre          \pagetitretrue

\footline={\ifpagetitre\the\baspagetitre\else
\ifodd\pageno\the\baspagedroite\else\the\baspagegauche\fi\fi
\global\pagetitrefalse}

\def\raggedbottom{\topskip 10pt plus 36pt\r@ggedbottomtrue}

\def\pointir{\unskip . --- \ignorespaces}

\def\Bigbreak{\vskip-\lastskip\bigbreak}
\def\Medbreak{\vskip-\lastskip\medbreak}

\def\ctexte#1\endctexte{%
  \hbox{$\vcenter{\halign{\hfill##\hfill\crcr#1\crcr}}$}}

\long\def\ctitre#1\endctitre{%
    \ifdim\lastskip<24pt\vskip-\lastskip\bigbreak\bigbreak\fi
  		\vbox{\parindent=0pt\leftskip=0pt plus 1fill
          \rightskip=\leftskip
          \parfillskip=0pt\bf#1\par}
    \bigskip\nobreak}

\long\def\section#1\endsection{%
\vskip 0pt plus 3\normalbaselineskip
\penalty-250
\vskip 0pt plus -3\normalbaselineskip
\Bigbreak
\message{[section \string: #1]}{\bf#1\unskip}\pointir}

\long\def\sectiona#1\endsection{%
\vskip 0pt plus 3\normalbaselineskip
\penalty-250
\vskip 0pt plus -3\normalbaselineskip
\Bigbreak
\message{[sectiona \string: #1]}%
{\bf#1}\medskip\nobreak}

\long\def\subsection#1\endsubsection{%
\Medbreak
{\it#1\unskip}\pointir}

\long\def\subsectiona#1\endsubsection{%
\Medbreak
{\it#1}\par\nobreak}

\let\+=\tabalign

\def\signature#1\endsignature{\vskip 15mm minus 5mm\rightline{\vtop{#1}}}

\mathcode`A="7041 \mathcode`B="7042 \mathcode`C="7043 \mathcode`D="7044
\mathcode`E="7045 \mathcode`F="7046 \mathcode`G="7047 \mathcode`H="7048
\mathcode`I="7049 \mathcode`J="704A \mathcode`K="704B \mathcode`L="704C
\mathcode`M="704D \mathcode`N="704E \mathcode`O="704F \mathcode`P="7050
\mathcode`Q="7051 \mathcode`R="7052 \mathcode`S="7053 \mathcode`T="7054
\mathcode`U="7055 \mathcode`V="7056 \mathcode`W="7057 \mathcode`X="7058
\mathcode`Y="7059 \mathcode`Z="705A

\def\spacedmath#1{\def\packedmath##1${\bgroup\mathsurround=0pt ##1\egroup$}%
\mathsurround#1 \everymath={\packedmath}\everydisplay={\mathsurround=0pt }}

\def\nospacedmath{\mathsurround=0pt \everymath={}\everydisplay={} }

\def\decale#1{\smallbreak\hskip 28pt\llap{#1}\kern 5pt}
\def\decaledecale#1{\smallbreak\hskip 34pt\llap{#1}\kern 5pt}
\def\puce{\smallbreak\hskip 6pt{$\scriptstyle\bullet$}\kern 5pt}

\def\displaylinesno#1{\displ@y\halign{
\hbox to\displaywidth{$\@lign\hfil\displaystyle##\hfil$}&
\llap{$##$}\crcr#1\crcr}}

\def\ldisplaylinesno#1{\displ@y\halign{
\hbox to\displaywidth{$\@lign\hfil\displaystyle##\hfil$}&
\kern-\displaywidth\rlap{$##$}\tabskip\displaywidth\crcr#1\crcr}}

\def\eqalign#1{\null\,\vcenter{\openup\jot\m@th\ialign{
\strut\hfil$\displaystyle{##}$&$\displaystyle{{}##}$\hfil
&&\quad\strut\hfil$\displaystyle{##}$&$\displaystyle{{}##}$\hfil
\crcr#1\crcr}}\,}

\def\system#1{\left\{\null\,\vcenter{\openup1\jot\m@th
\ialign{\strut$##$&\hfil$##$&$##$\hfil&&
        \enskip$##$\enskip&\hfil$##$&$##$\hfil\crcr#1\crcr}}\right.}

\let\@ldmessage=\message

\def\message#1{{\def\pc{\string\pc\space}%
                \def\'{\string'}\def\`{\string`}%
                \def\^{\string^}\def\"{\string"}%
                \@ldmessage{#1}}}

\def\diagram#1{\def\normalbaselines{\baselineskip=0pt
\lineskip=10pt\lineskiplimit=1pt}   \matrix{#1}}

\def\up#1{\raise 1ex\hbox{\smallf@nt#1}}

\def\qed{\raise -2pt\hbox{\vrule\vbox to 10pt{\hrule width 4pt
                 \vfill\hrule}\vrule}}

\def\cqfd{\unskip\penalty 500\quad\vrule height 4pt depth 0pt width
4pt\medbreak}

\def\virg{\raise .4ex\hbox{,}}   


\def\build#1_#2^#3{\mathrel{
\mathop{\kern 0pt#1}\limits_{#2}^{#3}}}

\def\boxit#1#2{%
\setbox1=\hbox{\kern#1{#2}\kern#1}%
\dimen1=\ht1 \advance\dimen1 by #1 \dimen2=\dp1 \advance\dimen2 by #1
\setbox1=\hbox{\vrule height\dimen1 depth\dimen2\box1\vrule}%
\setbox1=\vbox{\hrule\box1\hrule}%
\advance\dimen1 by .4pt \ht1=\dimen1
\advance\dimen2 by .4pt \dp1=\dimen2  \box1\relax}

\def\date{\the\day\ \ifcase\month\or janvier\or f\'evrier\or mars\or
avril\or mai\or juin\or juillet\or ao\^ut\or septembre\or octobre\or
novembre\or d\'ecembre\fi \ {\old \the\year}}

\def\crog{{\vrule height 2.57mm depth 0.85mm width 0.3mm}\kern -0.36mm
[}

\def\crod{]\kern -0.4mm{\vrule height 2.57mm depth 0.85mm
width 0.3 mm}}

\def\rond{\kern 1pt{\scriptstyle\circ}\kern 1pt}

\def\diagram#1{\def\normalbaselines{\baselineskip=0pt
\lineskip=10pt\lineskiplimit=1pt}   \matrix{#1}}

\def\hfl#1#2{\nospacedmath\smash{\mathop{\hbox to
12mm{\rightarrowfill}}\limits^{\scriptstyle#1}_{\scriptstyle#2}}}

\def\phfl#1#2{\nospacedmath\smash{\mathop{\hbox to
8mm{\rightarrowfill}}\limits^{\scriptstyle#1}_{\scriptstyle#2}}}

\def\vfl#1#2{\llap{$\scriptstyle#1$}\left\downarrow\vbox to
6mm{}\right.\rlap{$\scriptstyle#2$}}

\catcode`\@=12

\showboxbreadth=-1  \showboxdepth=-1
\magnification=1200

\baselineskip=14pt
\spacedmath{1.7pt}
\baspagegauche={\centerline{\tenbf\folio}}
\baspagedroite={\centerline{\tenbf\folio}}
\hautpagegauche={\hfil}
\hautpagedroite={\hfil}
\def\pa{\S\kern.15em}
\def\ra{\rightarrow}
\def\saut{\vskip 5mm plus 1mm minus 2mm}

\font\pc=cmcsc10 \rm
\def\ppav{principally polarized abelian variety }
\def\ppavs{principally polarized abelian varieties }
\def\pav{polarized abelian variety }

\def\char{{\rm char}(k)}
\def\Z{\hbox{\bf Z}}
\def\C{\hbox{\bf C}}
\def\End{\mathop{\rm End}\nolimits}
\def\Nm{\mathop{\rm Nm}\nolimits}
\def\Card{\mathop{\rm Card}\nolimits}
\def\rank{\mathop{\rm rank}\nolimits}

\def\Pic{\mathop{\rm Pic}\nolimits}
\def\Tr{\mathop{\rm Tr}\nolimits}
\def\Id{\hbox{\rm Id}}
\def\mod{\mathop{\rm mod}\nolimits}
\def\loc{{\it loc.cit.\/}}
\parskip=1.9mm
\def\ie{\hbox{i.e. }}
\def\ne{numerically equivalent }
\def\type{$(\delta_1|\cdots |\delta_n)$}
\def\dra{-\, -\ra}

\vskip 8mm
\ctitre
{\bf DEGREES OF CURVES IN ABELIAN VARIETIES}
\endctitre
\centerline{\pc Olivier Debarre
\footnote{(*)}{\rm Partially supported by the European Science
Project ``Geometry of Algebraic Varieties", Contract no. SCI-0398-C (A) and by
N.S.F. Grant DMS 92-03919.}
}
\medskip

{\pc R\'esum\'e} -- Le degr\'e d'une courbe $C$ contenue dans une vari\'et\'e
ab\'elienne polaris\'ee
$(X,\lambda)$ est l'entier $d=C\cdot\lambda$. Lorsque $C$ engendre $X$, on
obtient une minoration
de $d$ en fonction de $n$ et du degr\'e de la polarisation $\lambda$. Le plus
petit
degr\'e possible est $d=n$ et n'est atteint que pour une courbe lisse dans sa
jacobienne avec sa
polarisation principale canonique (Ran, Collino). On \'etudie les cas $d=n+1$
et $d=n+2$. Lorsque $X$
est simple, on montre de plus, en utilisant des r\'esultats de Smyth sur la
trace des entiers
alg\'ebriques totalement positifs, que si $d\le 1,7719\, n$, alors $C$ est
lisse et $X$ est isomorphe
\`a sa jacobienne. Nous obtenons aussi une borne sup\'erieure pour le genre
g\'eom\'etrique de $C$ en
fonction de son degr\'e.

{\pc Abstract} -- The degree of a curve $C$ in a \pav
$(X,\lambda)$ is the integer $d=C\cdot\lambda$. When $C$ generates $X$, we find
a lower bound on
 $d$ which depends on $n$ and the degree of the polarization $\lambda$. The
smallest possible
degree is $d=n$ and is obtained only for a smooth curve in its Jacobian with
its principal
polarization (Ran, Collino). The cases  $d=n+1$ and
$d=n+2$ are studied. Moreover, when $X$ is simple, it is shown, using results
of
Smyth on the trace of totally positive algebraic integers, that if $d\le
1.7719\, n$, then $C$ is
smooth and  $X$ is isomorphic to its Jacobian. We also get an upper bound on
the geometric genus of
$C$ in terms of its degree. \smallskip

{\bf 1. Introduction}

\ind Although curves in projective spaces have  attracted a lot of attention
for a long time, very
little is known in comparison about curves in abelian varieties. We try in this
article to
partially fill this gap.

\ind Let $(X,\lambda)$ be a \ppav of dimension $n$ defined over an
algebraically closed field $k$. The
degree of a curve $C$ contained in $X$ is $d=C\cdot \lambda$.

\ind The first question we are interested in is to find what numbers can be
degrees of curves $C$ in
$X$. When $C$ generates $X$, we prove that  $d\geq n(\lambda^n/n!)^{1/n}
\geq n$. It is known ([C], [R]) that $d=n$ if and only if $C$ is smooth
and $X$ is isomorphic to its Jacobian $JC$  with its  canonical principal
polarization. What
about the next cases? We get
 partial characterizations for $d=n+1$ and $d=n+2$, and we show (example 6.11)
that all degrees $\geq
n+2$ may actually occur (at least when $\char =0$). However, it seems necessary
to assume $X$ simple
to go further. We
prove, using results of Smyth ([S]), that {\it if $C$ is an irreducible curve
of degree $<2n$ if
$n\le 7$, and $\le 1.7719\, n$ if $n>7$, on a simple \ppav $X$ of dimension
$n$, then $C$
is smooth, has degree $2n-m$ for some divisor $m$ of $n$, the abelian variety
$X$ is isomorphic to
$JC$} (with a non-canonical principal polarization) {\it and $C$ is canonically
embedded in $X$.\/} We
conjecture this result to hold for any $n$ under the  assumption that $C$ has
degree $<2n$. This would
be a consequence of our conjecture 6.2, which holds for $n\le 7$: {\it the
trace of a totally positive
algebraic integer $\sigma$ of degree $n$ is at least $2n-1$ and  equality can
hold only if $\sigma$
has norm $1$}. Smooth curves of genus $n$ and degree $2n-1$ in their Jacobians
have
been constructed by Mestre for any $n$ in [Me].

\ind The second question is the Castelnuovo problem: bound the geometric genus
$p_g(C)$ of a curve
$C$ in a \pav $X$ of dimension $n$ in terms of its degree $d$. We prove, using
the original
Castelnuovo bound for curves in  projective spaces, the inequality
$p_g(C)<{(2d-1)^2\over 2(n-1)}$,
which is far from being sharp (better bounds are obtained for small degrees).
This in turn yields a
lower bound in $O(n^{3/2})$ on the degree of a curve in a {\it generic\/}
\ppav of dimension
$n$.

\ind Part of this work was done at the M.S.R.I. in Berkeley, and the author
thanks this
institution for its hospitality and support.\saut

{\bf 2. Endomorphisms and polarizations of abelian varieties}

\ind Let $X$ be an abelian variety of dimension $n$ defined over an
algebraically closed field $k$ and let $\End (X)$ be its ring of endomorphisms.
The degree $\deg (u)$
of an endomorphism $u$ is defined to be $0$ if $u$ is not surjective, and the
degree of $u$ as a map
otherwise. For any prime $l$ different from the characteristic of $k$, the Tate
module $T_l(X)$
is a free ${\bf Z}_l$--module of rank $2n$ ([Mu], p.171) and the $l$--adic
representation $\rho_l:\End (X)\ra \End (T_l(X))$ is injective. For any
endomorphism $u$ of $X$, the
characteristic polynomial of $\rho_l(u)$ has coefficients in ${\bf Z}$ and is
independent of $l$.
It is called the characteristic polynomial of $u$ and is denoted by $P_u$. It
satisfies
$P_u(t)=\deg (t\,\Id _X-u)$ for any integer $t$ ([Mu], theorem 4, p.180). The
opposite $\Tr (u)$ of
the coefficient of $t^{2n-1}$ is called the trace of $u$.

\ind The N\'eron-Severi group of $X$ is the group of algebraic equivalence
classes of invertible
sheaves on $X$. Any element $\mu$ of $NS(X)$ defines a morphism $\phi_\mu :X\ra
\Pic ^0(X)$
([Mu], p.60) whose scheme-theoretic kernel is denoted by $K(\mu )$. The
Riemann-Roch theorem gives
$\chi (X,\mu)=\mu^n/n!$, a number which will be called the degree of $\mu$. One
has $\deg
\phi_{\mu}=(\deg \mu)^2$ ([Mu], p.150). A polarization $\lambda$ on $X$ is the
algebraic equivalence
class of an ample invertible sheaf on $X$; it is said to be {\it separable\/}
if its degree is prime
to $\char$. In that case, $\phi_{\lambda}$ is a separable isogeny and its
kernel  is
isomorphic to a group $(\Z /\delta_1\Z )^2\times \cdots \times (\Z /\delta_n\Z
)^2$, where
$\delta_1|\cdots |\delta_n$ and $\delta_1\cdots \delta_n=\deg (\lambda)$. We
will say that $\lambda$
is of type $(\delta_1|\cdots |\delta_n)$.

\ind We will need the following result:
\medskip
{\pc Theorem} 2.1.(Kempf, Mumford, Ramanujan)-- {\it Let $X$ be an abelian
variety of dimension $n$,
and let $\lambda$ and $\mu$ be two elements of $NS(X)$. Assume that $\lambda$
is a polarization.
Then:}
{\parindent=2cm\item{\rm (i)}{\it The roots of the polynomial
$P(t)=(t\lambda-\mu)^n$ are all real.
\item {\rm (ii)} If $\mu$ is a polarization, the roots of $P$ are all positive.
\item {\rm (iii)} If $P$ has no negative roots and $r$ positive roots, there
exist a \pav $(X',\mu')$
of dimension $r$ and a surjective morphism $f:X\ra X'$ with connected kernel
such that
$\mu=f^*\mu'$.\par}}
\smallskip
{\bf Proof.} The first point is part of [MK], theorem 2, p.98. The second point
follows
from the same theorem and the fact that if   $M$ is an ample line bundle on $X$
with class $\mu$,
one has $H^i(X,M)=0$ for $i>0$ ([Mu], \S\kern.15em 16). For the last point, the
same theorem from
[MK] yields that the neutral component $K$ of the group $K(\mu)$ has dimension
$n-r$. The restriction
of $M$ to $K$ is algebraically equivalent to $0$ ({\it loc.cit.\/}, lemma 1,
p.95) hence, since the
restriction $\Pic ^0(X)\ra \Pic^0(K)$ is surjective, there exists a line bundle
$N$ on $X$
algebraically equivalent to $0$ such that the restriction of $M\otimes N$ to
$K$ is trivial. It
follows from theorem 1, p.95 of {\it loc.cit.\/} that $M\otimes N$ is the
pull-back of a line bundle
on $X'=X/K$.\cqfd

(2.2)\hskip .4cm Suppose now that $\theta$ is a principal polarization on $X$,
\ie a polarization
of degree $1$. It defines a morphism of ${\bf Z}$--modules $\beta_{\theta}:
NS(X)\ra \End (X)$
by the formula $\beta_{\theta}(\mu)= \phi_{\theta}^{-1}\circ \phi_{\mu}$. Its
image consists of all
endomorphisms invariant under the {\it Rosati involution\/}, which sends an
endomorphism $u$ to
$\phi_{\theta}^{-1}\circ\Pic ^0(u)\circ\phi_{\theta}$ ([Mu], (3) p.190).
Moreover, one has,
for any integer $t$:
$$\Bigl( {(t\theta-\mu)^n\over n!}\Bigr)^2=\deg
(t\phi_{\theta}-\phi_{\mu})=\deg \bigr( t\,\Id
_X-\beta_{\theta}(\mu)\bigr) =P_{\beta_{\theta}(\mu)}(t)\ .$$

(2.3)\hskip .4cm Let $(X,\lambda)$ be a polarized abelian variety. For $0<r\leq
n$, we set
$ \lambda^r_{\min}=$\break $\lambda^r/(r!\,\delta_1\cdots \delta_r)$.
If $k={\bf C}$, the class of $\lambda^r_{\min}$ is  minimal (\ie
non-divisible)  in $H^{2r}(X,{\bf Z})$. If $k$ is any algebraically closed
field, and if
$l$ is a prime number different from the characteristic of $k$, the group
$H^1_{\hbox{\sevenrm\'et}}(X,{\bf Z}_l)$ is a free ${\bf Z}_l$--module of rank
$n$ ([Mi], theorem
15.1) and the algebras $H^*_{\hbox{\sevenrm\'et}}(X,{\bf Z}_l)$ with its
cup-product structure and
$\wedge ^* H^1_{\hbox{\sevenrm\'et}}(X,{\bf Z}_l)$ with its wedge-product
structure, are isomorphic
([Mi], remark 15.4). In particular, the class $[\lambda]_l$ in
$H^2_{\hbox{\sevenrm\'et}}(X,{\bf
Z}_l)$ of the polarization $\lambda$ can be viewed as an alternating form on a
free ${\bf
Z}_l$--module, and as such has elementary divisors. If $\lambda$ is {\it
separable\/},
$(X,\lambda)$ lifts in characteristic $0$ to a \pav of the same type
$(\delta_1|\cdots |\delta_n)$.
The elementary divisors of $[\lambda]_l$ are then the maximal powers of $l$
that divide
$\delta_1,\ldots ,\delta_n$. Since intersection
corresponds to cup-product in \'etale cohomology, the class of
$\lambda^r_{\min}$ is in $H^{2r}_{\hbox{\sevenrm\'et}}(X,{\bf
Z}_l)$ and is not divisible by $l$.

\ind Throughout this article, all schemes we consider will be defined over  an
algebraically closed field $k$. We will denote numerical equivalence by $\sim$.
If $C$ is a
smooth curve, $JC$ will be its Jacobian and $\theta_C$  its canonical principal
polarization.\saut

{\bf 3. Curves and endomorphisms}

\ind We summarize here some results from [Ma] and [Mo]. Let $C$ be a curve on a
\pav $(X,\lambda)$
and let $D$ be an effective divisor that represents $\lambda$. Morikawa proves
that the following
diagram, where $d$ is the degree of $C$ and $S$ is the sum morphism, defines an
endomorphism
$\alpha(C,\lambda)$ of $X$ which is independent on the choice of $D$:
\vskip -8mm
$$\diagram{\alpha(C,\lambda):&X&\dra&C^{(d)}&\hfl{S}{}&X&\hfl{\rm
translation}{}&X\cr
{}&x&\longmapsto&(D+x)\cap C&{}&{}&{}&{}\cr}$$
\vskip -4mm
(3.1)\hskip .4cm Let $N$ be the normalization of $C$. The morphism
${\iota}:N\ra X$
factorizes through a morphism $p:JN\ra X$. Set $q=\iota^*\circ
\phi_{\lambda}:X\ra JN$; Matsusaka
proves that $ \alpha(C,\lambda)=p\circ q$ ([Ma], lemma 3).

(3.2)\hskip .4cm He also proves (\loc , theorem 2) that
$\alpha(C,\lambda)=\alpha(C',\lambda)$ if
and only if $C\sim C'$. Since $\alpha(\lambda^{n-1},\lambda)=(\lambda^n / n)\Id
_X$, it follows
that:
$$\alpha(C,\lambda)=m\,\Id _X \Longleftrightarrow C\sim {m \over (n-1)!\deg
\lambda}\lambda^{n-1}\ .$$

\ind If moreover $\lambda$ is separable of type $(\delta_1|\cdots |\delta_n)$
and if $l$ is a prime
distinct from $\char$, the discussion of (2.3) yields that there exists a class
$\epsilon$ in
$H^2_{\hbox{\sevenrm\'et}}(X,{\bf Z}_l)$ such that $\lambda^{n-1}\cdot
\epsilon$ is
$(n-1)!\delta_1\cdots\delta_{n-1}$ times a generator of
$H^{2n}_{\hbox{\sevenrm\'et}}(X,{\bf
Z}_l)$. It follows that $c=m / \delta_n$ must be in ${\bf
Z}_l$. But $\delta_n$ is prime to $\char$, hence $c$ is an integer and $C\sim
c\lambda^{n-1}_{min}$.

\ind Let $\theta_N$ be the canonical principal polarization on $JN$. One has:
$$\phi_{q^*\theta_N}=\Pic^0(q)\circ \phi_{\theta_N}\circ q=\phi_{\lambda}\circ
p\circ
\phi^{-1}_{\theta_N}\circ\iota^*\circ
\phi_{\lambda}=\phi_{\lambda}\circ\alpha(C,\lambda)\ .\leqno
(3.3)$$
Similarly:
$$
\phi_{p^*\lambda}=\Pic ^0(p)\circ \phi_{\lambda}\circ p=\phi_{\theta_N}\circ
q\circ\phi_{\lambda}\circ\phi_{\lambda}^{-1}\circ p= \phi_{\theta_N}\circ
q\circ p\ .$$
Note that, if $g$ is the genus of $N$, one has $C\cdot \lambda=N\cdot
p^*\lambda=\theta_N^{g-1}/(g-1)!\cdot p^*\lambda$. In particular, $-2(C\cdot
\lambda)$ is the coefficient
of $t^{2g-1}$ in the polynomial:$$
\deg (t\,\theta_N-p^*\lambda)^2=\deg
(t\,\phi_{\theta_N}-\phi_{p^*\lambda})=\deg (t\,\Id
_{JN}-q\circ p)\ . $$
Since $\Tr (q\circ p)=\Tr (p\circ q)=\Tr\bigl( \alpha(C,\lambda)\bigr)$,
the following equality, originally proved by Matsusaka ([Ma], corollary, p.8),
holds:$$ \Tr\bigl(
\alpha(C,\lambda)\bigr) = 2(C\cdot \lambda)\ .\leqno (3.4) $$

(3.5)\hskip .4cm If the N\'eron-Severi group of $X$ has rank $1$ (this holds
for a generic
principally polarized $X$ by [M], theorem 6.5, hence for a generic $X$ with any
polarization by
[Mu], corollary 1, p.234), and ample
generator $l'$, we can write $q^*\theta_N=rl'$ and $l=sl'$ with $r$ and $s$
integers. We get
$r\phi_{l'}=s\phi_{l'}\circ \alpha (C,l)$ hence $\alpha (C,l)(sx)=rx$ for all
$x$ in $X$. By taking
degrees, one sees that $s$ divides $r$ and $\alpha (C,l)=(r/ s)\Id _X$. By
(3.2), the curve $C$
is \ne to a rational multiple of $\lambda^{n-1}$ and its degree is a multiple
of $n$. If $l'$ is
separable of type \type, the curve $C$ is \ne to an integral multiple of
$\lambda^{n-1}_{\min }$ and
its degree is a multiple of $n\delta _n$.
\medskip
{\pc Lemma} 3.6.-- {\it Let $C$ be an irreducible curve that
generates a \pav $(X,\lambda)$ of dimension $n$. Then, the polynomial
$P_{\alpha(C,\lambda)}$ is the
square of a polynomial whose roots are all real and positive.}

\smallskip  {\bf Proof.} Let $\alpha=\alpha(C,\lambda)$. By (3.3), one has
$\phi_{q^*\theta_N}=\phi_{\lambda}\circ \alpha$, hence, for any integer $t$:
$$
\eqalign{P_{\alpha}(t)\ \deg \phi_{\lambda}&=\deg (t\,\Id _X-\alpha)\ \deg
\phi_{\lambda}\cr
{}&=\deg (t\,\phi_{\lambda}-\phi_{\lambda}\circ\alpha)\cr
{}&=\deg (t\,\phi_{\lambda}-\phi_{q^*\theta_N})\cr
{}&=\deg (\phi_{t\lambda-q^*\theta_N})= \bigl[ {1\over n!}\
(t\lambda-q^*\theta_N)^n\bigr]^2\ .\cr}$$
The lemma then follows from theorem 2.1.\cqfd

\ind We end this section with a proof of Matsusaka's celebrated criterion:
\medskip
{\pc Theorem} 3.7.(Matsusaka)-- {\it Let $C$ be an irreducible curve in a \pav
$(X,\lambda)$ of dimension $n$. Assume that $\alpha(C,\lambda)=\Id _X$. Then
$C$
is smooth and $(X,\lambda)$ is isomorphic to $(JC,\theta_C)$.}

\smallskip  {\bf Proof.} Let $N$ be the normalization of $C$. The morphism
$\alpha(C,\lambda)$ is
the identity and factors as:
$$
X\dra N^{(n)}\longrightarrow W_n(N)\longrightarrow JN\longrightarrow X\ .
$$
It follows that $\dim JN=g(N)\geq n$. Moreover, the image of $X$ in $JN$ has
dimension $n$, hence is
the entire $W_n(N)$, which is therefore an abelian variety. This is possible
only if $g(N)\leq n$.
Hence $N$ has genus $n$. It follows that the morphism $q:X\ra JN$ is an
isogeny, which is in fact an
isomorphism since $p\circ q=\alpha(C,\lambda)=\Id _X$. By (3.3), the
polarizations
$q^*\theta_N$ and $\lambda$ are equal, hence $q$ induces an isomorphism of the
polarizations.\cqfd

\saut
{\bf 4. Degrees of curves}

\ind Let $C$ be a curve that generates a \pav $(X,\lambda)$ of dimension $n$.
We want to study its
degree $d=C\cdot \lambda$. First, by  description (3.1), the dimension of the
image  of $\alpha(C,\lambda)$ is the dimension of the abelian subvariety $<C>$
generated by $C$. This
and  the definition of $\alpha(C,\lambda)$ imply:
$$
C\cdot \lambda\geq n\ .
$$
It was proved by Ran  ([R]) for $k={\bf C}$ and by Collino ([C]) in general,
that if $C\cdot
\lambda=n$, the minimal value, then $C$ is smooth and $(X,\lambda)$ is
isomorphic to its Jacobian
$(JC,\theta_C)$. This suggests that there should be a better lower bound on the
degree that involves
the type of the polarization $\lambda$. The following proposition provides such
a bound.
\medskip
{\pc Proposition} 4.1.-- {\it Let $C$ be an irreducible curve that generates a
\pav
$(X,\lambda)$ of dimension $n$. Then:}
$$
C\cdot \lambda\geq n(\deg \lambda)^{1 \over n}\ .
$$
{\it If $\lambda$ is separable, there is equality if and only if $C$ is smooth
and $(X,\lambda)$
is isomorphic to $(JC,\delta\theta_C)$, for some integer $\delta$ prime to
$\char$.\/}
\smallskip
\ind Recall that by (3.5), the degree of any curve on a {\it generic\/} \pav
$(X,\lambda)$ is a multiple of $n$. When $\lambda$ is separable of type \type,
this degree  is even a
multiple of $n\delta_n$.

\smallskip  {\bf Proof of the proposition.} We know by lemma 3.6 that
$P_{\alpha(C,\lambda)}$ is the
square of a polynomial $Q$ whose roots $\beta_1,\ldots ,\beta_n$ are real and
positive. We have:
$$
\eqalign{C\cdot \lambda &={1\over 2} \Tr \alpha(C,\lambda)={1\over 2}(2\beta_1
+\cdots +2\beta_n)\cr
&\geq n\ (\beta_1\cdots\beta_n)^{{1\over n}}=n\ Q(0)^{{1\over n}}=n\
P_{\alpha(C,\lambda)}(0)^{{1\over
2n}}\cr
&=n\ \bigl( \deg \alpha(C,\lambda)\bigr) ^{{1\over 2n}}\geq n\ ( \deg
\phi_{\lambda})
^{{1\over 2n}}=n\ (\deg\lambda)^{{1\over n}}\ .}$$
This proves the inequality in the proposition. If there is equality,
$\beta_1,\ldots ,\beta_n$ must
be all equal to the same number $m$, which must be an integer since
$P_{\alpha(C,\lambda)}$ has
integral coefficients. It follows from the proof of lemma 3.6 that:
$$
\bigl[ {1\over n!}\ (t\lambda-q^*\theta_N)^n\bigr]^2=P_{\alpha}(t)\ \deg
\phi_{\lambda}=(t-m)^{2n}\deg \phi_{\lambda}\ .
$$
Theorem 2.1.(iii) yields $m\lambda=q^*\theta_N$. It follows from (3.3) that
$\alpha(C,\lambda)=m\ \Id
_X$.

\ind If $\lambda$ is separable of type \type, by (3.2), the number $c=m/
\delta_n$ is  an
integer and $C$ is numerically equivalent to $c\lambda^{n-1}_{\min }$. We get:
$$cn\delta_n=C\cdot\lambda=n(\deg\lambda)^{{1\over
n}}=n(\delta_1\cdots\delta_n)^{{1\over n}}\le
n\delta_n\ .$$ This implies $c=1$ and $\delta_1=\cdots =\delta_n=\delta$. But
then $\lambda$ is
$\delta$ times a principal polarization $\theta$ ([Mu], theorem 3, p.231) and
$C\sim\theta^{n-1}_{\min }$. The conclusion now follows from Matsusaka's
criterion 3.7.\cqfd
\medskip
{\pc Corollary} 4.2.(Ran, Collino)-- {\it Let $C$ be an irreducible curve that
generates a \pav
$(X,\lambda)$ of dimension $n$. Assume that $C\cdot \lambda=n$. Then $C$ is
smooth and $(X,\lambda)$ is
isomorphic to its Jacobian $(JC,\theta_C)$.}
\smallskip
{\bf Proof.} Although the converse of the proposition was proved only for
$\lambda$ separable, we
still get from its proof that $\alpha(C,\lambda)$ is the identity of $X$ and we
may then apply
Matsusaka's criterion 2.7. This is the same proof as Collino's.\cqfd
\ind More generally, if $C\cdot\lambda=\dim <C>$, the same reasoning can be
applied on $<C>$ with
the induced polarization to prove that $C$ is smooth and that $(X,\lambda)$ is
isomorphic to the
product of $(JC,\theta_C)$ with a polarized abelian variety.
 \medskip {\pc Corollary} 4.3.-- {\it Let $X$ be an abelian variety with a
separable
 polarization $\lambda$ of type \type. Let $C$ be an irreducible curve in $X$
and let $m$ be
the dimension of the abelian subvariety that it generates.
Then:\/} $$
C\cdot \lambda\geq m(\delta_1\cdots\delta_m)^{1 \over m}\ .
$$

{\bf Proof.} Apply the proposition on the abelian subvariety $Y$ generated by
$C$. All there
is to show is that the degree $Y\cdot \lambda^m/ m!$ of the restriction
$\lambda'$ of $\lambda$
to $Y$ is at least $\delta_1\cdots\delta_m$. We will prove that it is actually
{\it divisible\/} by
$\delta_1\cdots\delta_m$. When $k=\C $, this follows from the fact that the
class $\lambda^m_{\min}$
is integral. The following argument for the general case was kindly
communicated to the author by
Kempf. Let $\iota$ be the inclusion of $Y$ in $X$. Then $\phi_{\lambda'}=\Pic
^0(\iota)\circ\phi_{\lambda}\circ\iota$, hence $\deg (\lambda')^2$, which is
the order of the
kernel of $\phi_{\lambda'}$, is a multiple of the order of its subgroup
$K(\lambda)\cap Y$, hence a
fortiori a multiple of the order of its $(r,r)$ part $K'$. In other words,
since $K'\simeq (\Z
/\delta'_1\Z )^2\times \cdots \times (\Z /\delta'_m\Z )^2$ for some integers
$\delta'_1|\delta'_2|\ldots |\delta'_m$ prime to $\char $, it is enough to show
that
$\delta'_1\delta'_2\cdots \delta'_m$ is a multiple of $\delta_1\delta_2\cdots
\delta_m$.

\ind Let $l$ be a prime number distinct from $\char $ and let ${\bf F}_l$ be
the
field with $l$ elements. For any integer $s$, let $X_s$ be the kernel of the
multiplication by $l^s$
on $X$. Then $X_s/X_{s-1}$ is a ${\bf F}_l$--vector
space of dimension $2n$ of which $Y_s/Y_{s-1}$ is a subspace   of dimension
$2m$. Since
$K(\lambda)$ is isomorphic to $(\Z /\delta_1\Z )^2\times \cdots \times (\Z
/\delta_n\Z )^2$, the rank
over ${\bf F}_l$ of  $\bigl( K(\lambda)\cap X_s\bigr) /\bigl( K(\lambda)\cap
X_{s-1}\bigr) $
is twice the cardinality  of the set $\{ i\in \{1,\ldots
,n\}\, |\, l^s\ \hbox{divides}\ \,\delta_i \} $. The dimension
formula yields:
$$
\rank \bigl( K(\lambda)\cap Y_s\bigr) /\bigl( K(\lambda)\cap
Y_{s-1}\bigr) \ge 2\,\Card\,\{ i\, |\, l^s\ \hbox{divides}\ \,\delta_i \} - 2n
+ 2m\ .$$

But the rank of $\bigl( K(\lambda)\cap Y_s\bigr) /\bigl( K(\lambda)\cap
Y_{s-1}\bigr)
=(K'\cap X_s) /( K'\cap X_{s-1}) $ is also twice the
cardinality  of  $\{ i\in \{1,\ldots ,m\}\, |\, l^s\ \hbox{divides}\
\,\delta'_i\} $. It follows that:
$$\Card\, \{ i\in \{1,\ldots ,n\}\bigm| l^s\!\!\not| \,\,\delta_i\} \ge \Card\,
\{ i\in \{1,\ldots
,m\}\bigm| l^s\!\!\not| \,\,\delta'_i\}\ .$$

This   implies what we need.\cqfd

\medskip {\pc Corollary} 4.4.-- {\it Let $C$ be an irreducible curve that
generates a \ppav
$(X,\theta)$ of dimension $n$. Assume that $C$ is invariant by translation by
an element $\epsilon$
of $X$ of order $m$. Then $C\cdot\theta\geq n\,m^{1-{1\over n}}$.}

{\bf Proof.} Let $H$ be the subgroup scheme generated by $\epsilon$. The
abelian variety $X'=X/H$
has a polarization $\lambda$ of degree $m^{n-1}$ whose pull-back on $X$ is
$m\theta$ ([Mu],
corollary, p.231). If $C'$ is the image of $C$ in $X'$, the proposition yields
$C\cdot\theta=C'\cdot\lambda\geq n\,m^{1-{1\over n}}$.\cqfd

\ind Note that in the situation of corollary 4.4, if $(X,\theta)$ is a {\it
generic\/} \ppav of
dimension $n$, and $m$ is prime to $\char$, then $mn$ divides $C\cdot\theta$.
With the notation of
the proof above, this follows from the fact that any curve on $X'$ is
numerically equivalent to an
integral multiple of $\lambda^{n-1}_{\min }$ (see (3.5)). \saut

{\bf 5. Bounds on the genus}

\ind We keep the same setting: $C$ is a curve that generates a \pav
$(X,\lambda)$, its normalization
is $N$, and its degree is $d=C\cdot\lambda$. The composition:
$$
X\dra N^{(d)}\longrightarrow W_d(N)\longrightarrow JN
$$
is a morphism with finite kernel (since $\alpha(C,\lambda)$ is an isogeny),
hence $W_d(N)$ contains
an abelian variety of dimension $n$. We can apply the ideas of [AH] to get a
bound of Castelnuovo
type on the genus of $N$. Note that if $C$ does not generate $X$, the same
bound holds with $n$
replaced by the dimension of $<C>$.
\medskip
{\pc Theorem} 5.1.-- {\it Let $C$ be an irreducible curve that generates a
separably \pav
$(X,\lambda)$ of dimension $n>1$. Let $N$ be the normalization of $C$ and let
$d=C\cdot\lambda$.
Then:\/}
$$
g(N) < {(2d-1)^2\over 2(n-1)}\ .
$$
\ind The inequality in the second part of lemma 8 in [AH] would improve this
bound
when $\char =0$, but its proof is incorrect.    \smallskip

{\bf Proof.} Let $A$ be the image of $X$ in $W_d(N)$ and let $A_2$ be the image
of $A\times A$ in $W_{2d}(N)$
under the addition map. We want to show that the morphism associated with a
generic point of $A_2$ is generically injective on $N$. The linear systems
corresponding to points
of $A_2$ are of the form $|{\cal O}_N(2D_x)|$, where $x$ varies in $X$, where
$D$ is an effective
divisor that represents $\lambda$ and $D_x=D+x$. It is therefore enough to show
that the restriction
to $C-x$ of the morphism $\phi_{2D}$ associated with $|2D|$ is generically
injective for $x$ generic.
If not, for $x$ generic in $X$ and $a$ generic in $C-x$, there exists $b$ in
$C-x$ with $a\neq b$ and
$\phi_{2D}(a)=\phi_{2D}(b)$. The same holds for $a$ generic in $X$ and $x$
generic in $C-a$. Since
$\phi_{2D}$ is finite ([Mu], p.60), $b$ does not depend on $x$, hence
$C-a=C-b$. Since $C$
generates $X$, this implies that $\epsilon=a-b$ is torsion, hence does not
depend on $a$. Letting $a$
vary, we see that any divisor in $|2D|$ is invariant by translation by
$\epsilon$. The argument in
[Mu], p.164, yields a contradiction.

\ind It follows that the image of the morphism $N\ra {\bf P}^r$ that
corresponds to a generic
point in $A_2$ is a curve of degree a divisor $d'$ of $2d$, with normalization
$N$. Moreover,
one has $r\geq n$ ([AH], lemma 1). Castelnuovo's bound ([ACGH], p.116 and [B]
when $\char >0$)
then gives:
$$
g(N)\leq m(d'-1)-{m(m+1)(r-1)\over 2}\ ,
$$
where $m=\bigl[ {d'-1\over r-1}\bigr]$. Hence:
$$
\eqalign{
g(N) &\leq m\biggl( d'-1-{(m+1)(r-1)\over 2}\biggr)\cr
&< {d'-1\over r-1}\biggl( d'-1-{d'-1\over 2}\biggr)\cr
&\leq {(d'-1)^2\over 2(n-1)}\leq {(2d-1)^2\over 2(n-1)}\ .\cr}
$$
This finishes the proof of the theorem.\cqfd

\ind In particular, in a \ppav $(X,\lambda)$ of dimension $n$, any {\it
smooth\/} curve numerically
equivalent to $c\theta^{n-1}_{\min}$ has genus $<{(2cn-1)^2\over 2(n-1)}$. For
curves in {\it
generic\/} \ppavs of dimension $n$, I conjecture the stronger inequality
$g(C)\leq cn+(c-1)^2$.

\ind The theorem also gives a lower bound on the degree of any curve in a  {\it
generic\/} complex
\pav of dimension $n$, whose only merit is to go to infinity with $n$. A better
bound is obtained in
[D2]. \medskip
{\pc Corollary} 5.2.-- {\it Let $C$ be a curve in a generic complex \pav
$(X,\lambda)$ of dimension $n$ and let $c$ be the integer such that $C$ is
\ne to  $c\,\lambda^{n-1}_{\min}$. Then:\/}
$$
c> \sqrt{{n\over 8}}-{1\over 4}\ .
$$
\smallskip
{\bf Proof.} We may assume that $\lambda$ is a principal polarization and that
$n>12$. Let $N$ be the
normalization of $C$. Corollary 5.5 in [AP] yields $g(N)> 1+n(n+1)/ 4$, which,
combined
with the proposition, gives what we want.\cqfd

\ind We can get better bounds on the genus when $d/ n$ is small.

\medskip
{\pc Proposition} 5.3.-- {\it Let $C$ be an irreducible curve that generates a
complex  \pav
$(X,\lambda)$ of dimension $n$. Let
$N$ be the normalization of $C$ and let $d=C\cdot\lambda$. Then:}
{\parindent=2.5cm\item{\rm (i)}{\it If $d<2n$, then $g(N)\leq d$.
\item {\rm (ii)} If $d=2n$, then $g(N)< {3d\over 2}=3n$.
\item {\rm (iii)} If $d\leq 3n$, then $g(N)\leq 4d$.\item {\rm (iv)} If $d\leq
4n$, then $g(N)\leq
6d$.\par}}
\smallskip

{\bf Proof.} We keep the notation of the proof of  theorem  5.1. In particular,
$W_d(N)$
contains an abelian variety $A$ of dimension $n$. If $2n>d$, it follows from
proposition 3.3 of
[DF] that $g(N)\leq d$. Recall that we proved earlier that the morphisms that
correspond to
generic points in $A_2$ are {\it birational\/} onto their image. It follows
from corollary 3.6
of \loc\ that $g(N)<3d/ 2$ when $d=2n$. This proves (ii). We will do (iv) only,
(iii) being
analogous. First, we may assume that the embedding of $A$ in $W_d(N)$ satisfies
the minimality
assumptions made in [A1]. Let $A_k$ be the image of $A\times\cdots\times A$ in
$W_{kd}(N)$ under
the addition map and let $r_k$ be the maximum integer such that $A_k$ is
contained in
$W^{r_k}_{kd}(N)$. If $g(N)>6d$, we get, as in the proof of proposition 3.8 of
[DF], the inequalities
$r_6\geq 8n+2$ and $n\leq 6d-3r_6$. It follows that $d\geq (n+3r_6)/
6\geq(25n+6)/ 6>4n$. This proves
(iv).\cqfd

\ind The  inequality (ii) should be compared with the inequality $g(C)\leq
2n+1$ proved by Welters
in [W] when $\char =0$ for any irreducible curve $C$ \ne to
$2\theta^{n-1}_{\min}$ on a \ppav
$(X,\theta)$ of dimension $n$ (so that $C\cdot\theta=2n$). Equality is obtained
only with the Prym
construction.
\saut
{\bf 6. Curves of low degrees}

\ind Let $C$ be an irreducible curve that generates a \ppav
$(X,\theta)$ of dimension $n$. We keep the same notation: $N$ is the
normalization of $C$ and
$q:X\ra JN$ is the induced morphism. From (2.2), we get that the square of the
monic polynomial
$Q_C(T)=(T\theta-q^*\theta_N)^n/n!$ has integral coefficients (and is the
characteristic
polynomial of $\alpha(C,\theta)$). It follows that $Q_C$ itself has integral
coefficients, and we
get from theorem 2.1 and (3.4):
{\parindent=2cm\item{(i)}The roots of $Q_C$ are all real and positive.
\item{(ii)} The sum of the roots of $Q_C$ is the degree $d=C\cdot\theta$.
\item{(iii)} The product of the roots of $Q_C$ is the degree of the
polarization $q^*\theta_N$.\par}

\ind Smyth obtained in [S] a lower bound on the trace of a totally real
algebraic integer in terms
of its degree. His results can be partially summarized as follows.
\medskip
{\pc Theorem} 6.1.(Smyth)-- {\it Let $\sigma$ be a totally positive algebraic
integer of
degree $m$. Then $
\Tr (\sigma)>1.7719\,m$, unless $\sigma$ belongs to an explicit finite set, in
which
case $\Tr (\sigma)=2m-1$ and $\Nm (\sigma)=1$.} \smallskip
\ind It is tempting to conjecture:
\medskip
(6.2)\hskip .4mm {\pc Conjecture $C_m$} -- {\it Let $\sigma$ be a totally
positive algebraic
integer of degree $m$. Then $\Tr (\sigma)\geq 2m-1$. If there is equality, then
$\Nm (\sigma)=1$.}

(6.3)\hskip 4mm The inequality in the conjecture follows from Smyth's theorem
for $m\leq 8$ (and
holds also for $m=9$ according to further calculations). Smyth also worked out
a list of all totally
positive algebraic integers $\sigma$ for which $\Tr (\sigma)-\deg (\sigma)\leq
6$. It follows from
this list that {\it the full conjecture holds for $m\leq 7$.}

\ind There are infinitely many examples for which the conjectural bound is
obtained: if $M$ is an odd
prime, the algebraic integer $4\cos ^2(\pi/ 2M)$ is totally positive, has
degree $(M-1)/2$,
trace $M-2$ and norm $1$.
\medskip
{\pc Proposition} 6.4.-- {\it Let $C$ be an irreducible curve that generates a
 \ppav
$(X,\theta)$ of dimension $n$ and let $Q_C$ be the polynomial defined above.
Then, if $|Q_C(0)|=1$,
the curve $C$ is smooth, $X$ is isomorphic to its Jacobian and $C$ is
canonically embedded.\/}
\smallskip

{\bf Proof.} By fact (iii) above, the polarization $q^*\theta_N$ is principal.
The proposition then
follows from the following lemma.\cqfd
\medskip
{\pc Lemma} 6.5.-- {\it Let $(JN,\theta_N)$ be the Jacobian of a smooth curve,
let $X$ be a
non-zero abelian variety and let $q:X\ra JN$ be a morphism. Assume that
 $q^*\theta_N$ is a principal polarization. Then $q$ is an isomorphism.\/}
\smallskip
{\bf Proof.} Since
 $q^*\theta_N$ is a principal polarization, $q$ is an closed immersion. By
Mumford's proof of
Poincar\'e's complete reducibility theorem ([Mu], p.173), there exist another
abelian subvariety $Y$
of $JN$ and an isogeny $f:X\times Y\ra JN$ such that  $f^*\theta_N$ is the
product of the induced
polarizations on each factor. As in \loc , for any $k$--scheme $S$, the set
$(X\cap Y)(S)$ is
contained in $K(q^*\theta_N)(S)$, which is trivial. Hence $f$ is an isomorphism
of polarized
varieties. But a Jacobian with its canonical principal polarization cannot be a
product, hence $Y$ is
$0$ and $q$ is an isomorphism.\cqfd

\ind We now give a result on curves on {\it simple\/} abelian varieties. The
part that depends on the
validity of  conjecture 6.2  holds in particular for $n\leq 7$.
\medskip
{\pc Theorem} 6.6.-- {\it Let $C$ be an irreducible curve in a simple   \ppav
$(X,\theta)$ of dimension $n$. Assume that either $C\cdot\theta \le
1.7719\, n$, or that conjecture $C_m$ holds
for all divisors $m$ of $n$ and $C\cdot\theta <2n$. Then,  the curve $C$ is
smooth, $X$ is isomorphic
to its Jacobian and $C$ is canonically embedded.\/}
\smallskip
{\bf Proof.} Since $X$ is simple, the polynomial $P_{\alpha(C,\theta)}$, hence
also its ``square
root'' $Q_C$, are powers of an irreducible polynomial $R$ of degree some
divisor $m$ of $n$. If
the degree of $C$, which is equal to the sum of the roots of $Q_C$,   is $ \le
1.7719\, n$, the sum of the roots of $R$ is also $ \le
1.7719\, m$. It follows from theorem 6.1 that $|R(0)|=1$. On the other hand, if
$C\cdot\theta
<2n$, the sum of the roots of $R$ is also $<2m$, hence, since
$C_m$ is supposed to hold, we also have $|R(0)|=1$. The theorem then follows in
both cases from
proposition 6.4.\cqfd

\ind It follows from the proof of the theorem that $C$ has degree $2n-m$ for
some divisor $m$ of
$n$. In particular, for $n$ {\it prime\/}, either $C$ has degree $n$ and
$\theta$ is the canonical
principal polarization, or it has degree $2n-1$.

\ind If one wants curves of degree between $n$ and $2n$ in a simple abelian
variety $X$, and if one
believes in conjecture 6.2, $X$ needs to be {\it a Jacobian with real
multiplications\/} (in the
sense that the ring $\End (X)\otimes{\bf Q}$ contains a totally real number
field different from ${\bf
Q}$). Examples have been constructed in [Me] (see also [TTV]). More precisely,
for any integer
$M\geq 4$, Mestre constructs an explicit $2$--dimensional family of complex
hyperelliptic Jacobians
$JC$ of dimension $[ M/ 2] $ whose endomorphism rings contain a subring
isomorphic to
${\bf Z}[T]/G_M(T)$, where:
$$
G_M(T)=\prod_{0<k\leq [M/2]} \bigl( T-4\cos ^2{k\pi \over M}\bigr) \ ,
$$
whose elements are invariant under the Rosati involution. By (2.2), they
correspond to polarizations
on $JC$. Take $M$ {\it odd\/} and set $n=\dim (JC)=(M-1)/ 2$. Then, the
endomorphism of $X$ that
corresponds to $T$ gives rise to a {\it principal\/} polarization on $JC$, with
respect to which the
degree of $C$, canonically embedded, is $2n-1$. Therefore, for any $n\geq 2$,
we have examples
of {\it complex \ppavs of dimension $n$ that contain curves of degree
$2n-1$\/}. They are simple if
$2n+1$ is prime. For $n=2$, these examples are Humbert surfaces, which contain
curves of degree $3$
([vG], p.221).

\ind If the assumption $X$ simple is dropped,  much less  can be said. If $Q$
is a monic
polynomial with integral coefficients whose roots are all real, we will say
that a curve $C$ has
{\it real multiplications by $Q$\/} if there is an endomorphism of $JC$ whose
characteristic
polynomial (see \S\kern.15em 2) is $Q^2$. If $k={\bf C}$, this is the same as
asking that the
characteristic polynomial of the endomorphism acting on the space of
first-order differentials of $C$
be $Q$.
\medskip
{\pc Proposition} 6.7.-- {\it Let $C$ be an irreducible curve that generates a
 \ppav
$(X,\theta)$ of dimension $n$. Then,
if $C\cdot\theta =n+1$, the curve $C$ is smooth, $X$ is isomorphic to its
Jacobian and $C$ is
canonically embedded. Moreover, the curve $C$ has real multiplications by
$(T-1)^{n-2}(T^2-3T+1)$.}
\smallskip
{\bf Proof.} By theorem 6.1 and Smyth's list in [S], the polynomial $Q_C$ can
only be
$(T-1)^{n-1}(T-2)$ or $(T-1)^{n-2}(T^2-3T+1)$. By proposition 6.4, we need only
 exclude the first
case. By theorem 2.1, there exist a polarized elliptic curve $(X',\lambda')$
and a morphism $f':X\ra
X'$ such that $f'^*\lambda'=q^*\theta_N-\theta$. Similarly, there exist an
$(n-1)$--dimensional \pav
$(X'',\lambda'')$ and a morphism $f'':X\ra X''$ such that
$f''^*\lambda''=2\theta-q^*\theta_N$. The
isogeny $(f',f''):(X,\theta)\ra (X',\lambda')\times (X'',\lambda'')$ is a
morphism of polarized
abelian varieties. Since $\theta$ is principal, it is an isomorphism and
$\lambda'$ and $\lambda''$
are both principal polarizations. Then, $(X,q^*\theta_N)$ is isomorphic to
$(X',2\lambda')\times
(X'',\lambda'')$. In particular, the pull-back of $\theta_N$ by $X''\ra JN$ is
a principal
polarization. By lemma 6.5, this cannot occur.\cqfd

\ind In the next case where $\deg (C)=n+2$, the same techniques give partial
results.
\medskip
{\pc Proposition} 6.8.-- {\it Let $C$ be an irreducible curve that generates a
 \ppav
$(X,\theta)$ of dimension $n>2$. Assume that $\char \neq 2,3$. Then,
if $C\cdot\theta =n+2$, one of the following possibilities occurs:}
{\parindent=2cm\item{\rm (i)}{\it the curve $C$ is smooth of genus $n$, $X$ is
isomorphic to its Jacobian and $C$ is canonically embedded. Moreover, the curve
$C$ has real
multiplications by $(T-1)^{n-3}(T^3-5T^2+6T-1)$ or $(T-1)^{n-4}(T^2-3T+1)^2$.
\item {\rm (ii)} the curve $C$ is smooth of genus $n$ and bielliptic, \ie there
exists a morphism
 of degree $2$ from $C$ onto an elliptic curve $E$. The abelian variety $X$ is
the quotient
of $JC$ by an element of order $3$ that comes from $E$. \item {\rm (iii)} The
normalization $N$ of
$C$ has genus $n$ and real multiplications by\break $(T-1)^{n-2}(T^2-4T+2)$.
There is an isogeny
$JN\ra X$ of degree $2$, and either $C$ is smooth, or it has one node and $N$
is hyperelliptic. \item
{\rm (iv)} the curve $C$ is smooth and bielliptic of genus $n+1$, and has real
multiplications by
$T(T-1)^{n-2}(T^2-4T+2)$. The abelian variety $X$ is the ``Prym variety''
associated with the
bielliptic structure.\par}}
\medskip
{\pc Remark} 6.9.  Mestre's construction for $N=7$ gives examples of curves of
degree $5$ in
\ppavs of dimension $3$, which fit into case (i) of the proposition. Example
6.11 below show that case (ii) does occur. These are the only examples I know
of.

\smallskip {\bf Proof.} By theorem 5.1 and Smyth's
list in [S], the polynomial $Q_C$ can only be $(T-1)^{n-2}(T-2)^2$,
$(T-1)^{n-1}(T-3)$,
$(T-1)^{n-3}(T^3-5T^2+6T-1)$,\break $(T-1)^{n-2}(T^2-4T+2)$ or
$(T-1)^{n-4}(T^2-3T+1)^2$. If the
constant term is $\pm 1$, the same proof as above yields that we are in case
(i). The first
polynomial is excluded as in proposition 6.7 (use $n>2$).

\ind If $Q_C(T)=(T-1)^{n-1}(T-3)$, as in the proof of proposition 6.7, there
exist a
polarized elliptic curve $(X',\lambda')$ and a morphism $f':X\ra X'$ with
connected kernel $X''$ such
that $f'^*\lambda'=q^*\theta_N-\theta$ or equivalently
$q^*\theta_N=\theta+(\deg \lambda')[X'']$. The
identity ${1\over n!}(T\theta-q^*\theta_N)^n=(T-1)^{n-1}(T-3)$ yields $(\deg
\lambda')(\deg
\theta_{|X''})=2$. If $\deg \lambda'=2$, one gets a contradiction as in the
proof of proposition 6.7.
If $\lambda'$ is principal, one has $\deg \bigl( (q^*\theta_N)
_{|X''}\bigr)=2$. We use the following
result. \medskip
{\pc Lemma} 6.10.-- {\it Let $(JN,\theta_N)$ be the Jacobian of a smooth curve,
let $X$ be a
non-zero abelian variety and let $r:X\ra JN$ be a morphism with finite kernel.
Assume that
 $\deg (r^*\theta_N)$ is $\leq \dim (X)$ and prime to $\char$. Then $g(N)<\dim
(X)+\deg
(r^*\theta_N)$.\/} \smallskip {\bf Proof.} Let $K$ be the kernel of $r$ and let
$\iota :X/K\ra JN$
be the induced embedding. By Poincar\'e's complete reducibility theorem ([Mu],
p.173), there exist
an abelian subvariety $X'$ of $JN$ and an isogeny $f:X/K\times X'\ra JN$ such
that the pull-back
$f^*\theta_N$ is the product of the induced polarizations. Note that $\deg
(\iota^*\theta_N)$
divides $\deg (r^*\theta_N)$. In particular, under our assumptions, the
polarization
$\iota^*\theta_N$ is {\it separable\/} and has a non-empty base locus $F$, of
dimension $\geq \dim
(X)-\deg (r^*\theta_N)$. If $\Theta$ is a theta divisor for $JN$, it follows
from the equation of
$f^*\Theta$ given in [D1], proposition 9.1, that $f(F\times X')$ is contained
in $\Theta$. Lemma 5.1
from [DF] (which is valid in any characteristic) then yields $\dim (F\times
X')+\dim (X')\leq
g(N)-1$, from which the lemma follows.\cqfd

\ind Since $\char \neq 2$, it follows from the lemma applied to the inclusion
$X''\ra JN$ that
$g(N)=n$ hence that the morphism $q: X\ra JN$ is an isogeny of degree $3$. It
is not difficult to see
(using for example [D1] \S\kern.15em 9) that since $\char\neq 2,3$, there is a
commutative diagram
of separable isogenies:
\vskip -.7cm plus 1mm minus 2mm
$$
\diagram{X''\times X'&\hfl{3:1}{}&X''\times E&\hfl{3:1}{}&X''\times X'\cr
\vfl{4:1}{}&&\vfl{4:1}{}&&\vfl{4:1}{}\cr
X&\hfl{q}{}&JN&\hfl{p}{}&X\cr}
$$
\vskip -.5cm plus 1mm minus 2mm
where $E$ is the quotient of $X'$ by a subgroup of order $3$. The middle
vertical arrow  induces an
injection of  $E$ into $JN$ whose image has degree $2$ with respect to
$\theta_N$. By
duality, one gets a morphism $f:N\ra E$ of degree $2$. In this situation, one
checks that since
$n>2$, for any two points $x$ and $y$ of $N$, one cannot have $x-y\equiv f^*e$,
for $e\neq 0$ in $E$.
Thus $C$, image of $N$ in $X$ by $p$, is smooth.

\ind If $Q_C(T)=(T-1)^{n-2}(T^2-4T+2)$, the polarization $q^*\theta_N$ has
degree $2$. It follows
from lemma 6.10 that:
\parskip=0cm

\ind $\bullet$ either $g(N)=n$ and $C$ is the image of $N$ by an isogeny
$p:JN\ra X$ of degree $2$.
In particular, either $C$ is smooth or $N$ is hyperelliptic and $C$ is obtained
 by identifying two Weierstrass points of $N$ (so that, in a sense, $C$ is
bielliptic).

\ind $\bullet$ or $g(N)=n+1$ and $q$ is a closed immersion. The proof of lemma
6.10 yields an
elliptic curve $X'$
 in $JN$ and  an isogeny $f:X\times X'\ra JN$ of degree $4$. Moreover, $\deg
(\theta_N)_{|X'}=2$, hence the morphism $N\ra X'$ obtained by duality has
degree $2$. One checks as
above that $C$ is smooth. The abelian variety $X$ is the Prym variety
associated with the bielliptic
structure, \ie is isomorphic to the quotient $JN/X'$. It remains to prove the
statement about real
multiplications. With the notation of (2.2), we calculate the characteristic
polynomial of the
endomorphism $\beta_{\theta_N}(p^*\theta)$ of $JN$. If $t$ is any integer, one
has:
$$
\eqalign{
\deg \bigl( t\,\Id_{JN}-\beta_{\theta_N}(p^*\theta)\bigr) & = \deg
(t\,\theta_N-p^*\theta)^2 \cr
& = \bigl( {1 \over 4}\deg (t\, f^*\theta_N-f^*p^*\theta)\bigr) ^2\cr
& = \bigl( {1 \over 4}\deg (t\, (\theta_N)_{|X'})\deg (t\,
q^*\theta_N-q^*p^*\theta)\bigr) ^2\cr
 & = {t^2 \over 4}\deg \bigl( t\,\phi_{q^*\theta_N}-\phi_{q^*p^*\theta}\bigr)\
.\cr}
$$
Set $\alpha=\alpha(C,\theta)$. Using (3.1) and (3.3), we get:
$$
\eqalign{
\deg \bigl( t\,\Id_{JN}-\beta_{\theta_N}(p^*\theta)\bigr) & = {t^2 \over 4}\deg
\bigl(
t\,\phi_{\theta}\circ\alpha-\phi_{\alpha^*\theta}\bigr) \cr & = {t^2 \over
4}\deg \bigl( t\,\Id_{\Pic
^0(X)}-\Pic ^0(\alpha)\bigr) \,\deg (\phi_{\theta}\circ\alpha) \cr & = P_{\Pic
^0(\alpha)}(t)\, t^2=Q_C(t)^2\, t^2\ .\cr} $$ It follows that $N$ has real
multiplications by $T\,
Q_C(T)=T(T-1)^{n-2}(T^2-4T+2)$. This finishes the proof of the
proposition.\cqfd

{\pc Example} 6.11 Case (ii) of the proposition does occur as a particular case
of the following construction.
Let $C$ be a smooth curve of genus $n$ with a morphism of degree $r$ onto an
elliptic curve $E$.
Assume that $r$ is prime to $\char$ and that the induced morphism $E\ra JC$ is
a closed immersion.
Let $s$ be an integer prime to $\char$ and congruent to $1$ modulo $r$, and let
$q:JC\ra X$ be the
quotient  by a cyclic subgroup $H$
 of order $s$ of $E$. There exist an abelian variety $Y$ of dimension
$n-1$ with a polarization $\lambda_Y$ of type $(1|\cdots |1|r)$ and an isogeny
$f:E\times Y \ra
JC$ with kernel isomorphic to $({\bf Z}/r{\bf Z})^2$, such that
$f^*\theta_C=pr_1^*(r\lambda_E)\otimes pr_2^*\lambda_Y$, where $\lambda_E$ is
the polarization
on $E$ defined by a point. The isogeny $f$  induces an isogeny $g:E/H\times Y
\ra X$ and, because $s\equiv 1
\ (\mod\ r)$, one checks that there exists a  principal polarization $\theta$
on $X$ such that
$g^*\theta=pr_1^*(r\lambda_{E/H})\otimes pr_2^*\lambda_Y$. It follows that
$f^*q^*\theta=pr_1^*(rs\lambda_E)\otimes pr_2^*\lambda_Y$. We claim that {\it
the degree of
the curve $q(C)$  on $X$  with respect to the principal polarization $\theta$
is\/} $n+s-1$. In fact, one has:
$$
f^*\theta_C^{n-1}/(n-1)!\ \sim\
r\lambda_E\ (pr_2^*\lambda_Y)^{n-2}/(n-2)!+(pr_2^*\lambda_Y)^{n-1}/(n-1)!$$
hence
$$\eqalign{
f^*\theta_C^{n-1}/(n-1)!\cdot f^*q^*\theta &
=rs\deg\lambda_Y+r(n-1)\deg\lambda_Y\cr
& =r^2(s+n-1)\ .\cr}
$$
It follows that $C\cdot q^*\theta=n+s-1$, which proves our claim.

\ind When $\char =0$, this construction yields examples of curves of degree
$n+t$ in
Jacobians of dimension $n$, for any $n\geq 2$ and $t\geq 2$.

\saut  \centerline {\pc References}

\medskip

\hangindent=1cm
[A1]	D. Abramovich. -- Subvarieties of abelian varieties and of Jacobians of
curves. Ph.D. Thesis, Harvard University, 1991.

\hangindent=1cm
[A2]	D. Abramovich. -- Addendum to ``Curves and abelian
varieties on  $W_d(C)$''. Unpublished.

\hangindent=1cm
[AH]	D. Abramovich, J. Harris. -- Curves and abelian
varieties on  $W_d(C)$, {\it Comp. Math.} 78 (1991), 227--238.

\hangindent=1cm
[AP] 	A. Alzati, G.-P. Pirola. -- On abelian subvarieties generated by
symmetric
correspondences, {\it Math. Zeit.} 205 (1990), 333--342

\hangindent=1cm
[ACGH]	E. Arbarello, M. Cornalba, P. Griffiths, J. Harris. -- {\it Geometry Of
Algebraic Curves,
I}. Grundlehren 267, Springer Verlag, 1985.

\hangindent=1cm
[B]	E. Ballico. -- On singular curves in the case of positive characteristic,
{\it Math.
Nachr.} 141 (1989), 267--273.

\hangindent=1cm
[C]	A. Collino. -- A new proof of the Ran-Matsusaka criterion for Jacobians,
{\it Proc.
Amer. Math. Soc.} 92 (1984), 329--331.

\hangindent=1cm
[D1]	O. Debarre. -- Sur les vari\'et\'es ab\'eliennes dont le
diviseur th\^eta est singulier en codimension $3$, {\it Duke Math. J.} 56
(1988), 221--273.

\hangindent=1cm
[D2]	O. Debarre. -- Curves In Abelian Varieties. To appear.

\hangindent=1cm
[DF] O. Debarre, R. Fahlaoui. -- Abelian Varieties In  $W_d^r(C)$  And  Points
Of Bounded Degrees On Algebraic Curves. To appear in {\it Comp. Math.}

\hangindent=1cm
[Ma]	T. Matsusaka. -- On a
characterization of a Jacobian Variety, {\it Mem. Coll. Sc. Kyoto}, Ser. A,  23
(1959),
1--19.

\hangindent=1cm
[Me]	J.-F. Mestre. -- Familles de courbes hyperelliptiques \`a multiplications
r\'eelles. In {\it Arithmetic Algebraic Geometry}, Progress in Mathematics 89,
Birkh\"auser, 1991.

\hangindent=1cm
[Mi]	J. Milne. -- Abelian Varieties. In {\it Arithmetic Geometry}, edited
by G. Cornell and J. Silverman, Springer Verlag, 1986.

\hangindent=1cm
[M]	S. Mori. -- The endomorphism ring of some abelian varieties, {\it Japan J.
of Math.} 1 (1976),
109--130.

\hangindent=1cm
[Mo]	H. Morikawa. -- Cycles and Endomorphisms of Abelian Varieties, {\it Nagoya
Math.} J.7 (1954),
95--102.

\hangindent=1cm
[Mu] D. Mumford. -- {\it Abelian Varieties}. Oxford University Press, 1974.

\hangindent=1cm
[MK] 	D. Mumford. -- Varieties
Defined by Quadratic Equations, with an appendix by G. Kempf. In {\it Questions
On Algebraic
Varieties}, C.I.M.E., Varenna, 1970.

\hangindent=1cm
[R] Z. Ran. -- On
Subvarieties of Abelian Varieties, {\it Inv. Math.} 62 (1981), 459--479.

\hangindent=1cm
[S]	C. Smyth. -- Totally Positive Algebraic Integers of Small Trace, {\it Ann.
Inst. Fourier} 33
(1984), 1--28.

\hangindent=1cm
[TTV]	W. Trautz,
J. Top, A. Verberkmoes. -- Explicit hyperelliptic curves with real
multiplication and
permutation polynomials, {\it Can. J. Math.} 43 (1991), 1055--1064.

\hangindent=1cm
[vG]	G. van der Geer. -- {\it
Hilbert Modular Surfaces}. Ergebnisse der Math. und ihrer Grenz. 16, Springer
Verlag,
1988.

\hangindent=1cm
[W]	G. Welters. -- Curves with twice the minimal class on principally polarized
abelian varieties, {\it Proc. Kon. Ned. Akad. van Wetenschappen, Indagationes
Math.} 49
(1987), 87--109.
\vskip 1.2cm

O. {\pc Debarre}, Universit\'e Paris-Sud, Math\'ematique, B\^atiment 425, 91405
Orsay Cedex, France,
and Department of Mathematics, The University of Iowa, Iowa City, IA 52242,
U.S.A.
\medskip
A.M.S. Classification:
14K05, 14H40.\bye